\documentclass[10pt]{article}
\usepackage{amssymb}
\usepackage{amsmath,amsfonts,amssymb,amsthm}
\newtheorem{theorem}{Theorem}[section]

\newtheorem{lemma}[theorem]{Lemma}

\begin{document}

\title{Flag-transitive $4$-designs and $PSL(2,q)$ groups}

\author{Huili Dong\footnote{
Corresponding author:011141@htu.edu.cn }\\
{\small\it  College of Mathematics and Information Science, Henan Normal University,}\\
{\small\it Xinxiang, Henan 453007, P. R. China }\\
}
\date{}

\maketitle

\noindent{\bf Abstract}\quad    This paper considers   flag-transitive  $4$-$(q+1,k,\lambda)$ designs  with    $\lambda\geq5$  and $q+1>k>4$. Let
the automorphism  group of a design $\cal D$ be a  simple group  $G=PSL(2,q)$. Depend on  the  fact that the setwise stabilizer $G_B$ must be one of twelve   kinds of subgroups,   up
to isomorphism we get the following  two results.
(i) If $10\geq \lambda \geq 5$, then    except  $(G,G_x,G_B,k,\lambda)=(PSL(2,761),{E_{761}}\rtimes {C_{380}},S_4,24,7)$ or
$(PSL(2,512),{E_{512}}\rtimes {C_{511}},{D_{18}},18,8)$ undecided,  $\cal D$
 is   a  $4$-$(24,8,5)$, $4$-$(9,8,5)$,  $4$-$(8,6,6)$, $4$-$(10,9,6)$, $4$-$(9,6,10)$,     $4$-$(9,7,10)$, $4$-$(12,11,8)$ or  $4$-$(14,13,10)$   design with
 $G_B=D_8$, ${E_8}\rtimes {C_7}$, $D_6$, ${E_9}\rtimes {C_4}$, $PSL(2,2)$,  $D_{14}$, ${E_{11}}\rtimes {C_{5}}$ 
  or  ${E_{13}}\rtimes {C_6}$
  respectively.
 (ii)   If $\lambda>10$,  ${G_B}=A_4$, $S_4$,  $A_5$, $PGL(2,q_0)$($g>1$ even) or $PSL(2,q_0)$, where ${q_0}^g=q$, then
there is no such
  design.


\medskip \noindent{\bf MSC(2000) Subject Classification} 05B05, 05B25, 20B25

\section{Introduction}

A  $t$-$(v,k,\lambda)$ design  with a simple group $G=PSL(2,q)$ as an automorphism group  is  interesting. Many  classification results about  flag-transitive designs have  been  got.
 In 1986,  $2$-$(v,k,1)$ designs are classified by Delandtsheer\cite{Delandtsheer86}. In 2016,
 $2$-$(v,k,\lambda)$  symmetric designs  with the sole of $G$  equal to $PSL(2,q)$    are studyed  by Alavi, Bayat and Daneshkhah\cite{Alavi16}. In 2018,   all
nonsymmetric  $2$-designs with
$(r,\lambda)=1$   are determined by  Zhan and Zhou\cite{Zhan18}. In 2019, $2$-designs with $k=4$   are considered by   Zhan, Ding and Bai\cite{Zhan19}.
 For  $3$-designs, there are also  many results given in\cite{Cusack95, Keranen03,Keranen04,Liu12}. After   Dai and Li\cite{Dai17,Dai18} classified
  flag-transitive  $4$-designs with $\lambda=3$ or 4, we continue this work, and
get the following  two theorems:

\begin{theorem}  Let   $\cal D$  be a  flag-transitive $4$-$(q+1,k,\lambda)$ design with $10\geq \lambda \geq 5$ and   $q+1>k>4$, $G=PSL(2,q)$ be an automorphism simple group. Then
  up
to isomorphism except $(G,G_x,G_B,k,\lambda)=(PSL(2,761),{E_{761}}\rtimes {C_{380}},S_4,24,7)$ or
$(PSL(2,512),{E_{512}}\rtimes {C_{511}},{D_{18}},18,8)$ undecided,   $\cal D$
 is   a   $4$-$(24,8,5)$, $4$-$(9,8,5)$,  $4$-$(8,6,6)$, $4$-$(10,9,6)$, $4$-$(9,6,10)$,     $4$-$(9,7,10)$, $4$-$(12,11,8)$ or   $4$-$(14,13,10)$   design with
 $G_B=D_8$, ${E_8}\rtimes {C_7}$, $D_6$, ${E_9}\rtimes {C_4}$, $PSL(2,2)$,  $D_{14}$, ${E_{11}}\rtimes {C_{5}}$
  or  ${E_{13}}\rtimes {C_6}$
  respectively.
 
\end{theorem}

\begin{theorem}  Let $\cal D$ be  a flag-transitive  $4$-$(q+1,k,\lambda)$ design with $\lambda>10$ and  $q+1>k>4$, $G=PSL(2,q)$ be a simple  automorphism group of   $\cal D$, ${G_B}=A_4$, $S_4$,  $A_5$, $PGL(2,q_0)$($g>1$ even) or $PSL(2,q_0)$, where ${q_0}^g=q$. Then
there is no such
  design.
\end{theorem}

\section{Preliminaries}

We give some useful results.  The notations $n$,  $G_B$
and $G_{xB}$ denote $\gcd(2,q-1)$,  the setwise stabilizer    and $G_{x}\cap G_{B}$  respectively, where  $x\in B$.

\begin{lemma}  Let $\cal D$ be a  $4$-$(q+1,k,\lambda)$    design with  an  automorphism group $G=PSL(2,q)$ flag-transitively  acting on $\cal D$. Then

  \noindent(i) $\lambda (q-2)=\frac{k(k-1)(k-2)(k-3)}{n|G_{B}|}$;

 \noindent(ii) $q=\frac{(k-1)(k-2)(k-3)}{\lambda n|G_{xB}|}+2$;

 \noindent(iii) $k\mid (\lambda n(q-2)|G_{xB}|+6)$;

 \noindent(iv) $k\mid\gcd(\frac{q(q^2-1)}{n},\lambda n(q-2)|G_{xB}|+6)$.

\end{lemma}

{\bf Proof.}\, Since $G$ is flag-transitive, then by  \cite{Huber09}   $G$ is  $2$-transitive, therefore, $G$ is also    block-transitive and  point-primitive.
 From   $|G|=|G_{B}||B^G|=|G_{xB}||(x,B)^G|$,   we get  $b=\frac{|G|}{|G_{B}|}=\frac{|G|}{|G_{xB}|k}$. By \cite{CRC},
 $\frac{\lambda q(q+1)(q-1)(q-2)}{k(k-1)(k-2)(k-3)}=\frac{q(q^2-1)}{n|G_{B}|}=\frac{q(q^2-1)}{nk|G_{xB}|}$.
Obviously (i) and (ii) hold.  Since $(k-1)(k-2)(k-3)=(k^2-6k+11)k-6$, $k\mid|G_B|$ and  $|G_B|\mid|G|$,   (iii) and (iv) hold.

\bigskip

\section{Proof of two theorems}

 First we   assume that $\cal D$  is  a     $4$-$(q+1,k,\lambda)$   design with   $\lambda\geq5$ and  $q+1>k>4$,  $G=PSL(2,q)$  is
a flag-transitive  automorphism  group of $\cal D$, where $q=p^f\geq4$ and $p$ is a prime.
Since $\cal D$  is   flag-transitive, we get that   $k$ is one orbit length of $G_B$ acting on $\cal D$.
  By \cite{Dai18,Huber07}, up to conjugation  $G_B$ must be one of twelve   kinds of subgroups, then  several  results can be got.

\begin{lemma}      $G_B\neq A_4$, $S_4$ or $A_5$ except $(G,G_x,G_B,k,\lambda)=(PSL(2,761),$
${E_{761}}\rtimes {C_{380}},S_4,24,7)$ undecided.
\end{lemma}
{\bf Proof.}\, Assume that  $G_B= A_4$, $S_4$ or $A_5$, then $|G_{B}|=12$, $24$ or 60. Since $k\mid |G_{B}|$, we have $(|G_{B}|,k)=(12,12)$, $(12,6)$,
$(24,24)$, $(24,12)$, $(24,8)$, $(24,6)$, $(60,60)$, $(60,30)$, $(60,20)$, $(60,15)$, $(60,12)$,   $(60,10)$,  $(60,6)$ or $(60,5)$. Note that  $q>k-1$, $q$ is a power of a prime and
 $\lambda\geq5$,
by Lemma 2.1(i) and (ii), we get that
$(|G_{B}|,k,q,\lambda)=(12,12,32,33)$, $(12,12,13,45)$, $(12,12,17,33)$, $(12,12,47,11)$, $(12,12,101,5)$, $(12,6,8,5)$,
$(24,24,25,231)$, $(24,24,71,77)$, $(24,24,79,69)$, $(24,24,163,33)$,
 $(24,24,233,23)$, $(24,24,761,$
 $7)$, $(24,8,16,5)$, $(24,8,9,5)$,$(60,60,61,1653)$, $(60,60,89,1121)$, $(60,60,179,$
 $551)$, $(60,60,1123,87)$, $(60,30,128,87)$, $(60,30,31,189)$, $(60,30,89,63)$,
$(60,$
$30,263,21)$, $(60,30,191,29)$, $(60,20,53,19)$, $(60,20,59,$
$17)$, $(60,15,16,39)$, $(60,15,23,13)$, $(60,15,41,7)$, $(60,12,13,9)$ or
$(60,10,16,$
$6)$. By  the  lengths of orbits   listed in \cite{Dai18,Huber07}, we  need consider
$(|G_{B}|,k,q,\lambda)=(12,12,17,33)$, $(12,12,47,11)$, $(12,12,101,5)$,
 $(24,24,71,77)$, $(24,24,79,69)$,  $(24,24,233,2$
 $3)$, $(24,24,761,7)$, $(60,60,89,1121)$, $(60,60,179,551)$ or  $(60,30,89,63)$.
  Using Magma\cite{Magma} we can rule out all cases except case 7.  The possible  parameters  are listed in Table 1, where the notation $m^n$  means that the degree $m$ appears $n$ times.
 For example,   in case 1,   the  group  PrimitiveGroup(18,1) denotes  the primitive permutation group $G=PSL(2,17)$
 acting on the set $\Omega=\{1,2,\ldots,18\}$.
And $G$ have two  conjugacy classes of subgroups of  order $12$ and
the  possible lengths of all orbits of $G_B$  acting on  $\Omega$ are both
    $6$, 12. There are 2 possible basic blocks  $B_0$ of  designs.
 However, the command    Design$<$4,18$|$D$>$ where ${\cal D}={B}^G$ returns both structures are not 4-designs, then  case 1  can be  ruled out.
For case 7, from $|G_x|=\frac{|G|}{v}=289180$, we know $G_x={E_{761}}\rtimes {C_{380}}$. There are 64 possible basic blocks  $B_0$ of  designs.
 However, $b=\frac{|G|}{|G_{B}|}=9181465$  is too large to  use Magma. This case is undecided.

\begin{center}
{\rm
Table 1:  Possible  parameters    with $G_B=A_4$, $S_4$ or $A_5$}\bigskip

\begin{tabular}{|cc c c c|}
\hline

 $Case$   &$G_B$&  $(v,k,\lambda)$ &position   &  the lengths of all orbits\\
  \hline

1   & $A_4$   & $(18,12,33)$& (18,1)&6,12; 6,12 \\

2    & $A_4$   & $(48,12,11)$& (48,1)&$12^4$; $12^4$; $12^4$; $12^4$; $12^4$\\

3    & $A_4$   & $(102,12,5)$& (102,2)& 6, $12^8$ \\

4    & $S_4$   & $(72,24,77)$& (72,1)& $24^3$; $24^3$; $24^3$ \\

5    & $S_4$   & $(80,24,69)$& (80,1)& 8, $24^3$; 8, $24^3$\\
6   & $S_4$   & $(234,24,23)$& (234,4)& 6, $12$, $24^9$ \\

7 & $S_4$   & $(762,24,7)$& (761,1)& 6, 12, $24^{31}$; 6, 12, $24^{31}$\\

8 & $A_5$   & $(90,60,1121)$& (90,1)& 30, 60; 30, 60\\

9 & $A_5$   & $(180,60,551)$& (180,1)&  $60^3$; $60^3$; $60^3$\\

10 & $A_5$   & $(90,60,63)$& (90,1)&  30, 60; 30, 60\\

\hline

\end{tabular}

\end{center}

\bigskip

\begin{lemma}     $G_B=PGL(2,q_0)$, where ${q_0}^g=q$ and  $g>1$ even.
\end{lemma}
{\bf Proof.}\,Assume that  $G_B=PGL(2,q_0)$, where ${q_0}^g=q$ and  $g>1$ even.  Then
$k=q_0+1$, $q_0(q_0-1)$ or  $q_0({q_0}^2-1)$.

If $k=q_0+1$, then by Lemma 2.1(i)
$n\lambda ({q_0}^g-2)={q_0}-2$. However, there is no solution of  this equation.

If $k=q_0(q_0-1)$, then $q_0\geq 3$  and
\begin{align*}
  {q_0}^5-2> n\lambda ({q_0}^g-2)&=({q_0}-2)({q_0}^2-{q_0}-1)({q_0}^2-{q_0}-3).
    \end{align*}
   Thus $g=2$ or 4.
  Since
  \begin{align*}
  &({q_0}-2)({q_0}^2-{q_0}-1)({q_0}^2-{q_0}-3)\\
  &=({q_0}^2-2)[({q_0}-2)({q_0}^2-2{q_0}-1)]+({q_0}-2)\\
  &=({q_0}-4)({q_0}^4-2)+({q_0}^3+10{q_0}^2-3{q_0}-14),
  \end{align*}
  we get $({q_0}^2-2)\mid ({q_0}-2)$ or
     $({q_0}^4-2)\mid ({q_0}^3+10{q_0}^2-3{q_0}-14)$.  However,  this is impossible.

If $k=q_0({q_0}^2-1)$, then
\begin{align*}
 {q_0}^9-2> n\lambda ({q_0}^g-2)&=({q_0}^3-{q_0}-1)({q_0}^3-{q_0}-2)({q_0}^3-{q_0}-3),
   \end{align*}

   Thus $g=2$, 4, 6 or 8. Since
  \begin{align*}
&({q_0}^3-{q_0}-1)({q_0}^3-{q_0}-2)({q_0}^3-{q_0}-3)\\
&=({q_0}^2-2)({q_0}^7-{q_0}^5-6{q_0}^4+{q_0}^3+12{q_0}-6)+(13{q_0}-18)\\
&=({q_0}^4-2)({q_0}^5-3{q_0}^3-6{q_0}^2+5{q_0}+12)+(4{q_0}^3-18{q_0}^2-{q_0}+18)\\
 &=({q_0}^6-2)({q_0}^3-3{q_0}-6)+(3{q_0}^5+12{q_0}^4+12{q_0}^3-6{q_0}^2-17{q_0}-18)\\
&={q_0}({q_0}^8-2)+(-3{q_0}^7-6{q_0}^6+3{q_0}^5+12{q_0}^4+10{q_0}^3-6{q_0}^2-9{q_0}-6),
\end{align*}
 we get $({q_0}^2-2)\mid (13{q_0}-18)$, $({q_0}^4-2)\mid (4{q_0}^3-18{q_0}^2-{q_0}+18)$,   $({q_0}^6-2)\mid (3{q_0}^5+12{q_0}^4+12{q_0}^3-6{q_0}^2-17{q_0}-18)$ or   $({q_0}^8-2)\mid (3{q_0}^7+6{q_0}^6-3{q_0}^5-12{q_0}^4-10{q_0}^3+6{q_0}^2+9{q_0}+6)$.
   Therefore, $(g,{q_0})=(2,2)$, $(2,3)$.
 However,  this is contrary with $q>k-1$.

\begin{lemma}    Let  $G_B=PSL(2,q_0)$, where ${q_0}^g=q$. Then  $\cal D$
 is   a  $4$-$(9,6,10)$  design    with  $G_B=PSL(2,2)$, denoted by ${\cal D}_1$.
\end{lemma}
{\bf Proof.}\,Assume tht   $G_B=PSL(2,q_0)$, where ${q_0}^g=q$. Then $k=q_0+1$, $q_0(q_0-1)$ if $g$ is even or $\frac{{q_0}({q_0}^2-1)}{\gcd(2,{q_0}-1)}$.

If $k=q_0+1$, then    $\lambda({q_0}^g-2)=q_0-2$.  However,    there is no  such $q_0>3$ satisfying the  equation.

If $k=q_0(q_0-1)$, then $q_0\geq 3$  and
   ${q_0}^5-2>\lambda ({q_0}^g-2)=({q_0}-2)({q_0}^2-{q_0}-1)({q_0}^2-{q_0}-3)$.
   Thus $g=2$ or 4.
  The same as in Lemma 3.2,  this is impossible.

If $k=\frac{{q_0}({q_0}^2-1)}{\gcd(2,{q_0}-1)}$, then
   ${q_0}^9-2>{n^4}\lambda ({q_0}^g-2)=({q_0}^3-{q_0}-n)({q_0}^3-{q_0}-2n)({q_0}^3-{q_0}-3n)$. Note that  $q+1>k>4$,
   we get $8\geq g\geq 2$.

Since
\begin{align*}
&({q_0}^3-{q_0}-n)({q_0}^3-{q_0}-2n)({q_0}^3-{q_0}-3n)\\
&=({q_0}^2-2)[{q_0}^7-{q_0}^5-6n{q_0}^4+{q_0}^3+(11n^2+1){q_0}-6n]+[(11n^2+2){q_0}\\
&-6n^3-12n]\\
&=({q_0}^3-2)[{q_0}^6-3{q_0}^4-(6n-2){q_0}^3+3{q_0}^2+(12n-6){q_0}+(11n^2-12n\\
&+3)]+[-(6n-6){q_0}^2-(11n^2-24n+12){q_0}-6n^3+22n^2-24n+6]\\
&=({q_0}^4-2)({q_0}^5-3{q_0}^3-6n{q_0}^2+5{q_0}+12n)+[(11n^2-7){q_0}^3-18n{q_0}^2-\\
&(11n^2-10){q_0}-6n^3+24n]\\
&=({q_0}^5-2)({q_0}^4-3{q_0}^2-6n{q_0}+3)+
[(12n+2){q_0}^4+(11n^2-1){q_0}^3-(6n\\
&+6){q_0}^2-(11n^2+12n){q_0}-6n^3+6]\\
 &=({q_0}^6-2)({q_0}^3-3{q_0}-6n)+[3{q_0}^5+12n{q_0}^4+(11n^2+1){q_0}^3-6n{q_0}^2-\\
 &(11n^2+6){q_0}-6n^3-12n]\\
  &=({q_0}^7-2)({q_0}^2-3)+[-6n{q_0}^6+3{q_0}^5+12n{q_0}^4+(11n^2-1){q_0}^3-(6n-2)\\
  &{q_0}^2-11n^2{q_0}-6n^3-6]\\
&={q_0}({q_0}^8-2)+[-3{q_0}^7-6n{q_0}^6+3{q_0}^5+12n{q_0}^4+(11n^2-1){q_0}^3-6n{q_0}^2-\\
&(11n^2-2){q_0}-6n^3],
\end{align*}
 we get
 $(g,{q_0})=(3,2)$.  Thus $(q,k,\lambda)=(8,6,10)$.
The generators of   $G=PSL(2,8)$ acting on the set $\Omega=\{1,2,\ldots,9\}$ are  as follows:
$g_1=(1, 8)(2, 4)(3, 7)(5, 6)$, $g_2=(2, 7)(3, 6)(4, 5)(8, 9)$,   $g_3=(1, 2, 3, 4, 5, 6, 7)$.
 Then all orbits of $G_B$ acting  on  $\Omega$ are
$\Omega_{1}=\{ 4, 6, 8 \}$ and $\Omega_{2}=\{1, 2, 3, 5, 7, 9\}$.
Take
 $B=\Omega_{2}$ as a possible basic block of  a design, then it
 returns a  $4$-$(9,8,5)$    design  ${\cal D}_1$. Clearly, ${\cal D}_1$ is also flag-transitive.

\bigskip

Next we consider $G_B$ is one of the remaining four   kinds of subgroups and assume that $10\geq \lambda \geq 5$.
\bigskip

\begin{lemma} Let  $G_B={C_{c}}$ or $D_{2c}$, where $c\mid \frac{q\pm1}{n}$, then $(\lambda,|G_{B}|,c,k,q)=(5,8,4,8,23)$,  $(8,18,9,18,512)$, $(10,6,3,6,8)$,
$(7,8,8,8,17)$, $(7,8,4,8,17)$, $(10,14,7,7,8)$ or  $(6,6,3,6,7)$.

\end{lemma}
{\bf Proof.}\,Assume that  $G_B={C_{c}}$ or $D_{2c}$, where $c\mid \frac{q\pm1}{n}$.
Then
 $k=c$  or $2c$$(|G_B|=2c)$.

  Let  $c\mid \frac{q+1}{n}$,  then  by Lemma 2.1(iii) $c\mid\gcd(\lambda n(q-2)|G_{xB}|+6,\frac{q+1}{n})$, that is,
$c\mid\gcd(3\lambda n|G_{xB}|-6,\frac{q+1}{n})$. By Lemma 2.1(ii),  $(\lambda,|G_{B}|,c,k,q)=(5,8,4,8,23)$,  $(8,18,9,18,512)$ or $(10,6,3,6,8)$.

Let $c\mid \frac{q-1}{n}$.
Assume that  $(\lambda,n,|G_{B}|,k)=(6,1,2c,2c)$, then  $6(2^f-2)=(2c-1)(2c-2)(2c-3)=2c(4c^2-12c+11)-6$, therefore,   $3\frac{2^f-1}{c}=4c^2-12c+11$. From $3\mid (4c^2-12c+11)$, we have that
$c=3l+1$ or $3l+2$, where $l$ is a positive integer. If $c=3l+1$, then $2^{f}=l(6l+1)(6l-1)+2$, therefore, $8\mid (l(6l+1)(6l-1)+2)$, we get  $l\equiv 2 (\mod 8)$. Let
 $l=8m+2$, then  $2^{f-1}=(4m+1)(48m+13)(48m+11)+1$. Obviously,   $m\neq 0$,    therefore,  $2^{f-1}\geq5\cdot61\cdot59+1>2^{14}$.  Thus $2^{14}\mid ((4m+1)(48m+13)(48m+11)+1)$, however, this is impossible.
Assume that  $(\lambda,n,|G_{B}|,k)=(6,1,c,c)$, then
$6(2^f-2)=(c-1)(c-2)(c-3)=c(c^2-6c+11)-6$, therefore,   $6\frac{2^f-1}{c}=c^2-6c+11$. From $6\mid (c^2-6c+11)$, we have that
$c=6l+1$ or $6l+5$, where $l$ is an integer.
 If $c=6l+1$, then $2^{f-1}=l(6l-1)(3l-1)+1$, therefore, $4\mid (l(6l-1)(3l-1)+1)$, that is, $4\mid (18l^3-9l^2+l+1)$. However,  this is impossible.
If $c=6l+5$, then $2^{f-1}=(3l+2)(2l+1)(3l+1)+1$, therefore, $4\mid ((3l+2)(2l+1)(3l+1)+1)$, that is, $4\mid (18l^3+27l^2+13l+3)$. However,  this is also impossible.
Now we consider that     $(\lambda,n,|G_{B}|,k)\neq(6,1,2c,2c)$ and   $(6,1,c,c)$,
then $c\mid\gcd(\lambda n(q-2)|G_{xB}|+6,\frac{q-1}{n})$, that is,
$c\mid\gcd(\lambda n|G_{xB}|-6,\frac{q-1}{n})$. Thus
 $(\lambda,|G_{B}|,c,k,q)=(7,8,8,8,17)$, $(7,8,4,8,17)$, $(7,16,8,16,197)$,  $(10,14,7,7,8)$, $(6,6,3,6,7)$.
By employing the  command   Subgroups(G:OrderEqual:=b), we  rule out $(\lambda,|G_{B}|,c,k,q)=(7,16,8,16,197)$.

\begin{lemma} Let  $G_B=E_{q_0}$ or ${E_{q_0}}\rtimes {C_c}$, where $q_0\mid q$,  $c\mid (q_0-1)$  and $c\mid (q-1)$, then
$(\lambda,i,k,q)=(5,1,8,8)$,  $(6,2,9,9)$, $(8,2,11,11)$ or  $(10,2,13,13)$ where   $k=q_0$ and  $c=\frac{k-1}{i}$.
\end{lemma}
{\bf Proof.}\,Let  $G_B=E_{q_0}$ or ${E_{q_0}}\rtimes {C_c}$, where $q_0\mid q$,  $c\mid (q_0-1)$  and $c\mid (q-1)$.
Then
 $k=q_0$ or $c{q_0}$.

 Let $k=q_0$,       then  $k\mid\gcd(\lambda n(q-2)|G_{xB}|+6,q)$, that is, $k\mid\gcd(2\lambda n|G_{xB}|-6,q)$.
 If   $|G_{B}|=q_0$, then  $(\lambda,|G_{B}|,k,q)=(7,8,8,32)$.
 If   $|G_{B}|=c{q_0}$, then   $|G_{xB}|=c$, therefore,
 $k\leq 2\lambda nc-6$, so $\frac{k-1}{2\lambda n}<\frac{k+6}{2\lambda n}\leq c$.
  From $c\mid (k-1)$, we let $c=\frac{k-1}{i}$, where $i=1$, 2, $\cdots$, $2\lambda n-1$.
Then from    $k\mid\gcd(\frac{2\lambda n(k-1)}{i}-6, q)$,  we get   $k\mid \gcd(2\lambda n+6i,q)$.
All  possible parameters of   $(\lambda,i,k,q)$ are $(5,1,8,8)$, $(5,1,13,13)$, $(6,2,9,9)$,  $(6,3,7,7)$, $(7,1,17,17)$,
$(8,1,19,19)$, $(8,2,11,11)$, $(10,1,23,23)$ or  $(10,2,13,13)$.
By employing the  command   Subgroups(G:OrderEqual:=b), we  rule out    $(\lambda,i,k,q)=(5,1,13,13)$,   $(6,3,7,$
$7)$, $(7,1,17,17)$, 
$(8,1,19,19)$ and $(10,1,23,23)$.

 If $k=cq_0$,  then $|G_{xB}|=1$. From $k\mid(\lambda n(q-2)+6)$,
 we have that ${q_0}\mid\gcd(2\lambda n-6,q)$. Since  $c\mid (q_0-1)$. We get no possible parameters.

\begin{lemma}     Let  $G_B={C_{c}}$,  $D_{2c}$ or  $E_{q_0}$, ${E_{q_0}}\rtimes {C_c}$, where $c\mid \frac{q\pm1}{n}$ or  $q_0\mid q$,  $c\mid (q_0-1)$  and $c\mid (q-1)$. Then    except $(G,G_x,G_B,k,\lambda)=(PSL(2,512),{E_{512}}\rtimes {C_{511}},{D_{18}},18,8)$ undecided,   up
to isomorphism $\cal D$
 is   a  $4$-$(24,8,5)$, $4$-$(9,8,5)$,  $4$-$(8,6,6)$, $4$-$(10,9,6)$, $4$-$(9,6,10)$,     $4$-$(9,7,10)$, $4$-$(12,11,8)$   or   $4$-$(14,13,10)$   design with
 $G_B=D_8$, ${E_8}\rtimes {C_7}$, $D_6$, ${E_9}\rtimes {C_4}$, $PSL(2,2)$,  $D_{14}$, ${E_{11}}\rtimes {C_{5}}$    
  or  ${E_{13}}\rtimes {C_6}$
  respectively.

\end{lemma}
{\bf Proof.}\, We consider the cases appearing in Lemma 3.4 and 3.5 and   list the  results   in Table 2. The lengths of all orbits and the number of the designs are listed in
 column 5 and   column 6 respectively, where twice or five times denote such lengths of orbits appearing  twice or five times.

\begin{center}
{\rm
Table 2:  Possible values with $G_B={C_{c}}$, $D_{2c}$, $E_{q_0}$ or ${E_{q_0}}\rtimes {C_c}$}\bigskip

\begin{tabular}{|c  c c c cc|}
\hline

 $Case$    &$G_B$&  $(v,k,\lambda)$ &position   &   lengths& number\\
  \hline

1    & $D_8$   & $(24,8,5)$& (24,1)&$8^3$&2\\

2    & $D_{18}$   & $(513,18,8)$& (513,10)&9, ${18}^{28}$&?\\

3  & $D_6$   & $(9,6,10)$& (9,8)  &3,6&1\\

4  & $D_6$   & $(8,6,6)$&  (8,4)&2,6&1\\

5 & $D_8$   & $(18,8,7)$& (18,1)& 2, $4^2$, 8(twice)&0\\
6 & $C_8$   & $(18,8,7)$& (18,1)& $1^2$,$8^2$&0\\

7 & $D_{14}$   & $(9,7,10)$&(9,8) & 2,7&1\\

8 & $E_{8}$   & $(33,8,7)$&  (33,1)&  1, $8^4$(five times)&0\\

9 & ${E_8}\rtimes {C_7}$   & $(9,8,5)$& (9,8)&  1,8&1\\

10 & ${E_{9}}\rtimes {C_{4}}$   & $(10,9,6)$& (10,3)&  1,9&1\\

11 & ${E_{11}}\rtimes {C_{5}}$   & $(12,11,8)$& (12,1)&  1,11&1\\

12  & ${E_{13}}\rtimes {C_{6}}$   & $(14,13,10)$&(14,1) & 1,13 &1\\

\hline

\end{tabular}

\end{center}

For case 2,   from $|G_x|=\frac{|G|}{v}=261632$, we know $G_x={E_{761}}\rtimes {C_{380}}$. There are 28 possible basic blocks  $B_0$ of  designs.
 However, $b=\frac{|G|}{|G_{B}|}=7456512$ is too large to  use Magma.

For the remaining cases, the same as before,  we can deal with  them by  the same method.  For example,
   the  generators of the group $G=PSL(2,23)$ acting on the set $\Omega=\{1,2,\ldots,24\}$ are listed as follows:

 \noindent$g_1=
 (1, 8)(2, 10)(3, 16)(4, 24)(5, 15)(6, 21)(7, 23)(9, 20)(11, 22)(12,14)\\
 (13, 18)(17, 19)$, \\
 \noindent$g_2=(1, 19, 24, 6)(2, 13, 22, 20)(3, 12, 23, 5)(4, 17, 8, 21)(7, 14, 16,15)\\(9, 11, 18, 10)$, \\
 \noindent$g_3=(1, 24)(2, 22)(3, 23)(4, 8)(5, 12)(6, 19)(7, 16)(9, 18)(10, 11)(13,20)\\(14, 15)(17, 21)$.

 $G$ has only one  conjugacy class of subgroups of  order $8$ and
all the orbits of
 $G_B$  acting on the set $\Omega$  are listed as follows:

\noindent$\Omega_{1}=\{1, 4, 6, 8, 17, 19, 21, 24\}, \\
 \Omega_{2}=\{2, 9, 10, 11, 13, 18, 20, 22\}$, \\
 $\Omega_{3}=\{ 3, 5, 7, 12, 14, 15, 16, 23\}$.

\noindent  First take
 $B_0=\Omega_{1}$ as a possible basic block of  a design. However,
the structure is not 4-design.
 Now  take
 $B_0=\Omega_{2}$ or $\Omega_{3}$   as a possible basic block of  a design,
then
we  construct  two flag-transitive 4-$(24,8,5)$  designs  ${\cal D}_2$ and   ${\cal D}_3$.
  By\cite{Design}, it shows that both designs are isomorphic. It need to note that the design corresponding to case 3 is  isomorphic to  the design ${\cal D}_1$.

 This completes the proof of Theorem 1.1 and 1.2.

\noindent{\bf Acknowledgements}This work was supported by the
National Natural Science Foundation of China(No:11301158,11601132).

\end{document}